\documentclass[letterpaper, 10 pt, conference]{ieeeconf}  % Comment this line out
\IEEEoverridecommandlockouts                              % This command is only
\usepackage{graphicx} % for pdf, bitmapped graphics files
\usepackage{amsmath} % assumes amsmath package installed
\usepackage{amssymb}  % assumes amsmath package installed
\usepackage{url}

\newtheorem{thm}{Theorem}
\newtheorem{proposition}[thm]{Proposition}

\newtheorem{remark}{Remark}
\newtheorem{definition}{Definition}
\newtheorem{assumption}{Assumption}

\title{
Implementation of logical gates on infinite dimensional quantum oscillators
}

\author{\authorblockN{Nabile Boussa\"{i}d}
\authorblockA{Laboratoire de math\'ematiques\\
Universit\'e de Franche--Comt\'e\\
25030 Besan\c{c}on, France\\
{\tt\small Nabile.Boussaid@univ-fcomte.fr}}
\and
\authorblockN{Marco Caponigro}
\authorblockA{Institut \'Elie Cartan de Nancy and\\
INRIA Nancy Grand Est\\
54506 Vand{\oe}uvre, France\\
{\tt\small Marco.Caponigro@inria.fr}}
\and
\authorblockN{Thomas Chambrion}
\authorblockA{Institut \'Elie Cartan de Nancy and\\
INRIA Nancy Grand Est\\
54506 Vand{\oe}uvre, France\\
{\tt\small Thomas.Chambrion@inria.fr}}
}
\begin{document}

\maketitle
\thispagestyle{empty}
\pagestyle{empty}

 \begin{abstract}
In this paper  we study the error in the approximate simultaneous controllability of the bilinear Schr\"odinger equation.
We provide estimates based on a tracking algorithm for general bilinear quantum systems and on the study of the finite dimensional Galerkin approximations for a particular class of quantum systems, weakly-coupled systems. We then present two physical examples: the perturbed quantum harmonic oscillator and the infinite potential well.
\end{abstract}

\section{Introduction}
\subsection{Logical gates}
Quantum computation relies on the idea to store an information in the state of quantum system.
This state is described by the \emph{wave function}, that is, a point $\psi$ in the Hilbert sphere of  $L^2(\Omega,\mathbf{C})$, where $\Omega$ is a Riemannian manifold.

When submitted to an
excitation by an external field ({\it e.g.} a laser), the time evolution of the wave function is governed by the 
bilinear Schr\"odinger
equation %that reads
\begin{equation}\label{eq:blse}
i \frac{\partial \psi}{\partial t}=-\frac{1}{2}\Delta \psi +V(x)
\psi(x,t) +u(t) W(x) \psi(x,t),
\end{equation}
where $V, W:\Omega\rightarrow \mathbf{R}$ are real functions describing  respectively the
physical properties of the uncontrolled system and the external field, and $u:\mathbf{R}\rightarrow \mathbf{R}$ is a real function of the time representing the intensity of the latter.

When the manifold $\Omega$ is compact, the linear operator $i(\Delta/2 -V)$ admits a set of eigenstates $(\phi_n)_{n\in \mathbf{N}}$.%,\phi_2,\ldots,$.  
A \emph{logical gate}, or \emph{quantum gate}, is 
%nothing but 
a unitary transformation in $L^2(\Omega,\mathbf{C})$  for which some finite dimensional space of the form $\mathrm{span}\{\phi_1,\phi_2,\ldots,\phi_n\}$ is stable.
To build a given logical gate $\widehat{\Upsilon}$ from the system \eqref{eq:blse}, one has to find a control law $u$ such that the propagator $\Upsilon^u_T$ at a certain time $T$ of \eqref{eq:blse} satisfies $\Upsilon^u_T (\phi_j)=\widehat{\Upsilon} \phi_j$ for every $j=1,\ldots ,n$.

The main difficulty with this problem is that the space $L^2(\Omega,\mathbf{C})$ has infinite dimension.
% as soon as $\Omega$ has dimension  larger or equal to one. 
For the sake of simplicity, one often only considers the case where $\Omega$ is a finite union of points (or, equivalently, $L^2(\Omega,\mathbf{C})$ is finite dimensional). Nevertheless, most of the usual quantum systems evolves on non trivial manifolds $\Omega$. 
This papers deals with the effective implementation of some simple logical gates on models of quantum oscillators on 1-dimensional manifolds.

 \subsection{Quantum control}
 The problem of driving the solutions of \eqref{eq:blse} to a given target has been intensively studied in the past decades, both in the finite and infinite dimensional case. Many advances have been done in the infinite dimensional case, when there is only one source and one target. % (i.e., one is just interested in $\Upsilon^u_T \phi_1$). 
 The interested reader may refer, for instance and among many other references, to  \cite{camillo,beauchard-nersesyan} for the theoretical viewpoint and to \cite{Salomon} for numerical aspects. In particular, it was proved in \cite{BMS} that \emph{exact} controllability is impossible in general. This does not prevent to study \emph{approximate} controllability of \eqref{eq:blse}, that is to replace the condition $\Upsilon^u_T (\phi_j)=\widehat{\Upsilon} \phi_j$ by
 $\|\Upsilon^u_T (\phi_j)-\widehat{\Upsilon} \phi_j\|\leq \varepsilon$ for every $j=1,\ldots,n$.  To the best of our knowledge, there are only very few results of \emph{simultaneous} controllability
 %$\Upsilon^u_T \phi_1, \Upsilon^u_T \phi_2,\ldots,\Upsilon^u_T \phi_n$ 
 in the infinite dimensional case and the only available effective control techniques have been described in~\cite{Schrod2} and \cite{periodic}.

 Recently, we noticed in~\cite{feps} that a certain class of bilinear systems are precisely approached by their Galerkin approximations. Two important examples of these so-called \emph{weakly-coupled} systems are the quantum harmonic oscillator and the infinite potential well. The structure of weakly-coupled systems permits precise numerical simulations for the construction of quantum gates.

\subsection{Content of the paper}
The theoretical background is recalled in Section~\ref{sec:general}. Besides a quick survey on simultaneous control techniques for Equation \eqref{eq:blse}, we give precise definitions and approximation results for
weakly-coupled systems.  In Section~\ref{sec:oscillator} we apply these results to a perturbation of the quantum harmonic oscillator and we provide estimates for the error in the controllability in a suitable finite dimensional approximation. Similarly, in Section~\ref{sec:potentialwell} we study the infinite potential well.

\section{General theoretical results}\label{sec:general}

\subsection{Framework and notations}
We reformulate the problem \eqref{eq:blse} in a more abstract framework.  
%This will allow us  to treat examples slightly more general than \eqref{eq:blse}. 
In a separable Hilbert
space $H$
endowed with norm $\| \cdot \|$ and Hilbert product $\langle \cdot, \cdot
\rangle$, we consider the evolution problem
\begin{equation}\label{eq:main}
\frac{d \psi}{dt}=(A+ u B)\psi
\end{equation}
where $(A,B)$ satisfies Assumption~\ref{ass:ass}.

\begin{assumption}\label{ass:ass}
$(A,B)$ is a pair of linear operators such that
 \begin{enumerate}
  \item for every $u$ in $\mathbf{R}$, $A+uB$ is essentially
skew-adjoint;
\label{EgaliteDomaine}
  \item $A$ is skew-adjoint and has purely discrete spectrum $(-i
\lambda_k)_{k \in \mathbf{N}}$, associated to the Hilbert basis $(\phi_k)_{k\in \mathbf{N}}$ of eigenvectors of $A$.
 \end{enumerate}
\end{assumption}

From Assumption \ref{ass:ass}.\ref{EgaliteDomaine}, one deduces that, for every piecewise constant $u$, $u:t\mapsto =\sum_{j} u_j \chi_{(t_j,t_{j+1})}(t)$, with $0=t_0\leq t_1 \leq \ldots  \leq t_{N+1}$ and $u_0,\ldots,u_N$ in $\mathbf{R}$,
 the solution $t\mapsto \Upsilon^u_t\psi_0$ of \eqref{eq:main}  has the form
\begin{multline*}
 \Upsilon^u_t \psi_{0}=e^{(t-t_{j-1})(A + u_{j-1} B)}\circ\\\circ
e^{(t_{j-1}-t_{j-2})(A + u_{j-2} B)}\circ \cdots  \circ e^{t_{0}(A+
u_0  B)} \psi_{0},
\end{multline*}
for $t\in [t_{j-1},t_{j})$.

\begin{remark}
With extra regularity hypotheses, it is possible to define the propagator of \eqref{eq:main} for a larger class of controls. For instance, when $B$ is bounded, for every $t$, $\Upsilon_t:u\mapsto \Upsilon^u_t$ admits a unique continuous extension to $L^1_{\mathrm{loc}}(\mathbf{R},\mathbf{R})$.
\end{remark}

\subsection{Control results}
\begin{definition}
 Let $(A,B)$ satisfy Assumption~\ref{ass:ass}.
 A subset $S$ of
$\mathbf{N}^2$ \emph{couples} two levels $j,k$ in $\mathbf{N}$,
if
there exists a finite sequence $\big
((s^{1}_{1},s^{1}_{2}),\ldots,(s^{q}_{1},s^{q}_{2}) \big )$
in $S$ such that
\begin{description}
\item[$(i)$] $s^{1}_{1}=j$ and $s^{q}_{2}=k$;
\item[$(ii)$] $s^{j}_{2}=s^{j+1}_{1}$ for every $1 \leq j \leq q-1$;
\item[$(iii)$] $\langle  \phi_{s^{j}_{1}}, B \phi_{s^{j}_{2}}\rangle \neq 0$
for $1\leq j \leq
q$.
\end{description}

The subset $S$ is called a \emph{connectedness chain}   for $(A,B)$ if $S$
couples every pair of levels in $\mathbf{N}$.
A connectedness chain is said to be \emph{non-resonant} if for every $(s_1,s_2)$
in $S$, $|\lambda_{s_1}-\lambda_{s_2}|\neq |\lambda_{t_1}-\lambda_{t_2}|$ for
every $(t_{1},t_{2})$ in $\mathbf{N}^2\setminus\{(s_1,s_2),(s_2,s_1)\}$ such
that $\langle \phi_{t_{2}}, B \phi_{t_{1}}\rangle  \neq 0$.
\end{definition}

\begin{definition}
Let $(A,B)$ satisfy Assumption~\ref{ass:ass}. The system $(A,B)$ is
approximately simultaneously controllable if for every
$\widehat{\Upsilon} \in U(H)$ (unitary operators acting on $H$), $\psi_1,\ldots,\psi_n\in H$, and
$\varepsilon>0$, there  exists a piecewise constant function
$u_\varepsilon: [0,T_\varepsilon]\to \mathbf{R}$ such that
$$
\|\widehat{\Upsilon}\psi_j-\Upsilon_{T_\varepsilon}^{u_\varepsilon}
\psi_j\|<\varepsilon.
$$
for every $j=1,\ldots,n$.
\end{definition}

The following sufficient condition for approximate simultaneous controllability has been given in \cite{Schrod2}.

\begin{proposition}\label{prop:control}
 Let $(A,B)$  satisfy Assumption~\ref{ass:ass} and admit a non-resonant chain of connectedness.
Then $(A,B)$ is approximately simultaneously controllable.
\end{proposition}

\subsection{Weakly-coupled systems}

\begin{definition}
Let $k$ be a positive number and let  $(A,B)$ satisfy Assumption~\ref{ass:ass}.\ref{EgaliteDomaine}.
%be such that
% for every $u$ in $\mathbf{R}^{p}$, $A+uB$ is essentially
%skew-adjoint on $D(A)$, $(\lambda_k)_k$ tends to $+\infty$ and $i(A+uB)$ is bounded from
%below. 
Then $(A,B)$ is
\emph{$k$ weakly-coupled}
if for every $u\in \mathbf{R}$, $D(|A+uB|^{k/2})=D(|A|^{k/2})$ and
 there exists
a constant $C$ such that, for every $\psi$ in $D(|A|^k)$, $ |\Re \langle |A|^k
\psi,B_l\psi \rangle |\leq C |\langle |A|^k \psi, \psi \rangle|$.

%The \emph{coupling constant}  $c_k(A,B)$ of system $(A,B_1,\ldots,B_p)$ of
%order $k$ is the quantity
%$$
%\sup_{\psi\in D(|A|^{k})\setminus\{0\}} \sup_{1\leq l \leq p} \frac{ |\Re
%\langle |A|^k
%\psi,B_l\psi \rangle |}{|\langle |A|^k \psi, \psi \rangle|}.
%$$
\end{definition}

\begin{definition}
Let $N \in \mathbf{N}$. We define the projection $\pi_N: \psi\in H \mapsto \sum_{k\leq N}\langle \phi_k,\psi \rangle \phi_k$. The \emph{Galerkin approximation}  of \eqref{eq:main}
of order $N$ is the system in $H$
\begin{equation}\label{eq:sigma}
\dot x = (A^{(N)} + u B^{(N)}) x \tag{$\Sigma_{N}$}
\end{equation}
where $A^{(N)}=\pi_N A \pi_N$ and $B^{(N)}=\pi_N B \pi_N$ are the
\emph{compressions} of $A$ and $B$ (respectively).
\end{definition}

We denote by $X^{u}_{(N)}(t,s)$ the propagator of \eqref{eq:sigma}
associated with a piecewise constant functions $u$.

\begin{proposition}\label{prop:gga}
Let $k$ and $s$ be non-negative  numbers with
$0\leq s <k$. Let $(A,B)$ satisfy 
Assumption~\ref{ass:ass} and be $k$ weakly-coupled
%Assume that  for every $u$ in $\mathbf{R}$, $A+uB$ is essentially
%skew-adjoint on $D(A)$, that $i(A+uB)$ is bounded from
%below, and 
Assume that there
exist $d>0$, $0\leq r<k$ such that $\|B\psi \|\leq d \|\psi \|_{r/2}$ for
every $\psi$ in $D(|A|^{r/2})$.
Then
for every $\varepsilon > 0 $, $K\geq 0$, $n\in \mathbf{N}$, and
$(\psi_j)_{1\leq j \leq n}$ in $D(|A|^{k/2})^n$
there exists $N \in \mathbf{N}$
such that
for every piecewise constant function $u$
$$%\begin{equation}
\|u\|_{L^{1}} < K \implies \| \Upsilon^{u}_{t}(\psi_{j}) -
X^{u}_{(N)}(t,0)\pi_{N} \psi_{j}\|_{s/2} < \varepsilon,
$$%\end{equation}
for every $t \geq 0$ and $j=1,\ldots,n$.
\end{proposition}

\section{The perturbed quantum harmonic oscillator}\label{sec:oscillator}

\subsection{Physical model}
The quantum harmonic oscillator is one of the most studied quantum system. Schr\"{o}dinger equation reads
\begin{equation}\label{eq:harmonicoscillator}
i \frac{\partial \psi}{\partial t} (x,t) = -\frac{1}{2} \frac{\partial^{2} \psi}{\partial x^{2}} +
\left(\frac{1}{2} x^{2} - u(t)x \right)\psi(x,t)\,,
\end{equation}
where $x\in \Omega=\mathbf{R}$.
With the notations of \eqref{eq:main}, $A=-i(-\Delta + x^{2})/2$ and $B=i x$.

An Hilbert basis of $H$ made of eigenvectors of $A$ is given by the sequence of
the Hermite functions $(\phi_n)_{n \in \mathbf{N}}$, associated with the
sequence $( - i \lambda_n)_{n \in \mathbf{N}}$ of eigenvalues where
$\lambda_n=n-1/2$ for every $n$ in $\mathbf{N}$. In the basis $(\phi_n)_{n \in
\mathbf{N}}$, $B$ admits a tri-diagonal structure:
$$
\langle \phi_j,B\phi_k\rangle = \left \{\begin{array}{cl}
- i\sqrt{k-1} & \mbox{if } j=k-1\\
-i\sqrt{k} & \mbox{if } j=k+1\\
0 & \mbox{otherwise},
\end{array} \right.
$$
A chain of connectedness for this system is given by $S=\{(n,n+1)\,:\,n\in \mathbf{N}\}$. The chain $S$ is resonant indeed $|\lambda_{n+1}-\lambda_{n}| = 1$ for every $n$ in $\mathbf{N}$.
As a matter of fact, the system \eqref{eq:harmonicoscillator} is known to be non-controllable (see \cite{Rouchon,Illner}).

We consider a perturbation of this system.
Consider the inverse $A^{-1}$ of the operator $A$. The family $(\phi_n)_{n \in \mathbf{N}}$ is a family of eigenvectors for $A^{-1}$ associated with the eigenvalues  $( - i/ \lambda_n)_{n \in \mathbf{N}}$. For every $\eta \geq 0$ we set $A_{\eta} = A+\eta A^{-1}$. Since $A$ and $A^{-1}$ commute then  $(\phi_n)_{n \in \mathbf{N}}$ is a family of eigenvectors for $A_{\eta}$ associated with the eigenvalues
$( - i \lambda_n^{\eta})_{n \in \mathbf{N}}$ where $ \lambda_n^{\eta} = \lambda_n+\eta/\lambda_{n}$.
The set $S$ is a non-resonant chain of connectedness for system $(A_{\eta},B)$ for every $\eta >0$. Indeed
$\lambda^{\eta}_{n+1} - \lambda^{\eta}_{n} = 1-4 \frac{\eta}{4 n^{2} -1}$ and, clearly,
$\lambda^{\eta}_{n+1} - \lambda^{\eta}_{n}  = \lambda^{\eta}_{m+1} - \lambda^{\eta}_{m} $ if and only if
$n=m$.

By Proposition~\ref{prop:control} %\cite{Schrod2} 
the system $(A_{\eta},B)$ is approximately simultaneously controllable. Moreover by \cite[Theorem 2.13]{Schrod2} we have also an upper bound on the $L^{1}$-norm of the control independent of the error. For instance  we can steer approximately the first level $\phi_{1}$ to the second $\phi_{2}$  with a control law with $L^{1}$-norm smaller than
$5\pi/4$. Another consequence is that a quantum gate for $\phi_{1}, \phi_{2}$, and $\phi_{3}$ is approximately reachable, that is for every $\varepsilon>0$, there exists $t_{\varepsilon}>0$ and a piecewise constant function $u_{\varepsilon}$ such that $\|\Upsilon^{u_{\varepsilon}}_{t_{\varepsilon}}(\phi_j)-\phi_{\sigma(j)}\|<\varepsilon$ where $\sigma$ is the 3-cycle which exchanges $1,2$ and $3$. This can be achieved with  $\|u_{\varepsilon}\|_{L^1} \leq \pi/2( 1+\sqrt{2}/2)$.

\subsection{Estimates}
In the following, we only consider control of $L^1$-norm less than $K=\pi/2( 1+\sqrt{2}/2)$.
The particular tri-diagonal structure of system $(A_{\eta},B)$ is very useful for a priori estimates on the components of the propagator. Indeed if $\|u\|_{L^{1}} \leq K$, by~\cite[Remark 6]{feps}, we have that
\begin{equation}\label{eq:faketaylor}
|\langle \phi_{n+1},   \Upsilon^{u}_{t}(\phi_{j})\rangle| \leq \frac{(2K)^{n-2} }{(n-2)!} \sqrt{\frac{(2n-3)!}{(n-2)!}},
\end{equation}
for every $n$ in $\mathbf{N}$, $n\geq 3$ and $j=1,2,3$.

We use~\eqref{eq:faketaylor} to find estimates on the size $N$ of the Galerkin approximation whose existence is asserted by Proposition~\ref{prop:gga}.

First, let $N \geq j$  and notice that
\begin{eqnarray*}
\frac{d}{dt} \pi_{N} \Upsilon^{u}_{t}(\phi_{j})&=&
(A^{(N)} + u
B_l ^{(N)}) \pi_{N} \Upsilon^{u}_{t}(\phi_{j}) \\
&&\quad \quad \quad+
u (t) \pi_{N} B (\mathrm{Id} - \pi_{N}) \Upsilon^{u}_{t}(\phi_{j}).
\end{eqnarray*}
Hence, by variation of constants, for every $t \geq 0$,
\begin{eqnarray}\label{eq:variationofconstants}
\lefteqn{\pi_{N} \Upsilon^{u}_{t}(\phi_{j})= X^{u}_{(N)}(t,0)  \pi_{N}\phi_{j} } \nonumber \\
&&+
 \int_{0}^{t}
X^{u}_{(N)}(t,s) \pi_{N} B(\mathrm{Id} - \pi_{N}) \Upsilon^{u}_{s}(\phi_{j})
u(\tau)  d\tau.
\end{eqnarray}
Therefore, since $X^{u}_{(N)}(t,s)$ is unitary and for the tri-diagonal structure of the system we have, for $j=1,2,3$,
\begin{eqnarray*}
\lefteqn{\|\pi_{N} \Upsilon^{u}_{t}(\phi_{j}) - X^{u}_{(N)}(t,s) \phi_{j}\|}\\
&\leq &K \|\pi_{N} B(\mathrm{Id} - \pi_{N}) \Upsilon^{u}_{t}(\phi_{j})\| \\
& = &K| b_{N,N+1}||\langle \phi_{N+1}, \Upsilon^{u}_{t}(\phi_{j}) \rangle| \\
& = &K \sqrt{N}| \langle \phi_{N+1}, \Upsilon^{u}_{t}(\phi_{j}) \rangle| \\
& \leq &\frac{2^{N-1}K^{N-1} }{(N-2)!} \sqrt{\frac{(2N-3)!}{(N-3)!}}.
\end{eqnarray*}

Using $K=\pi/2( 1+\sqrt{2}/2) <2.69$, it is enough to consider a Galerkin approximation of size $N=420$ to get $\|\pi_{N} \Upsilon^{u}_{t}(\phi_{j}) - X^{u}_{(N)}(t,s) \phi_{j}\|\leq 10^{-4}$ for $j=1,2,3$.

\subsection{Numerical simulations}
For simulations, we choose $\eta=1$.
To induce the transition between levels $1$ and $2$, the control law we use is a piecewise constant $4\pi$ periodic function, taking value $1$ for $0\leq t < 5.\,10^{-3}$ and taking value 0 for $5.\,10^{-3} \leq t \leq 4\pi$. We apply this control for $314$ periods, that is during a time of $1256 \pi$ .

To induce the transition between levels 2 and 3, the control law we use is a piecewise constant $12\pi/5$ periodic function, taking value 1 for $0\leq t <5.\,10^{-3} $ and taking value $0$ for $5.\,10^{-3} \leq t < 12\pi/5$. We apply this control for 222 periods.

The simulations are done on a Galerkin approximation of size 420,  which garantees $\|\pi_{N} \Upsilon^{u}_{t}(\phi_{j}) - X^{u}_{(N)}(t,s) \phi_{j}\|\leq 10^{-4}$ for $j=1,2,3$.
At final time $T$, the resulting propagator is such that $|\langle \Upsilon^{u}_{T} \phi_1,\phi_3\rangle |>0.998$, $|\langle \Upsilon^{u}_{T} \phi_2,\phi_1\rangle |>0.999$ and $|\langle \Upsilon^{u}_{T} \phi_3,\phi_2\rangle |>0.999$. The time evolution of the moduli of the first coordinates of $\Upsilon^{u}_{t}(\phi_{j})$ for $j=1,2,3$ is depicted in Figures~1, 2, and 3.

All the computations were done using the free software NSP, see \cite{NSP}. The source code for the simulation is available at \cite{source2}.The total computation time is less than 4 minutes on a standard desktop computer.

\begin{figure}[tb!]
 \centering
 \includegraphics[width=8.4cm]{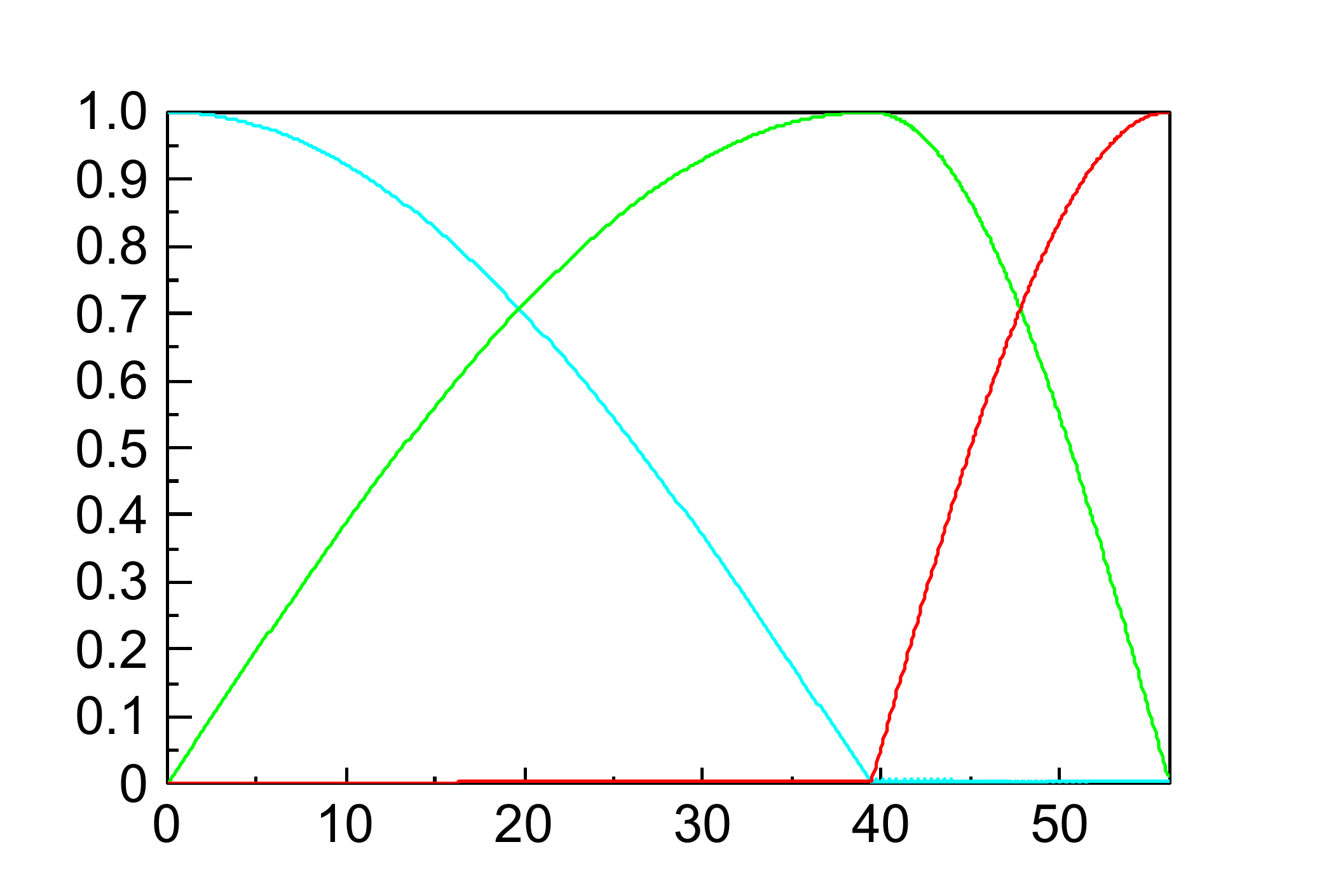}
 % HarmonicFromLevel1.pdf: 1276x696 pixel, 72dpi, 45.01x24.55 cm, bb=0 0 1276 696
 \caption{Time evolution of the moduli of the first three coordinates of $\Upsilon^u_t \phi_1$ in the case of the perturbed harmonic oscillator. First coordinate in blue, second coordinate in green, third coordinate in red. For the sake of readability, time scale is 1:100, total duration around 5500.}
 \label{fig:1}
\end{figure}

\begin{figure}[tb]
 \centering
 \includegraphics[width=8.4cm]{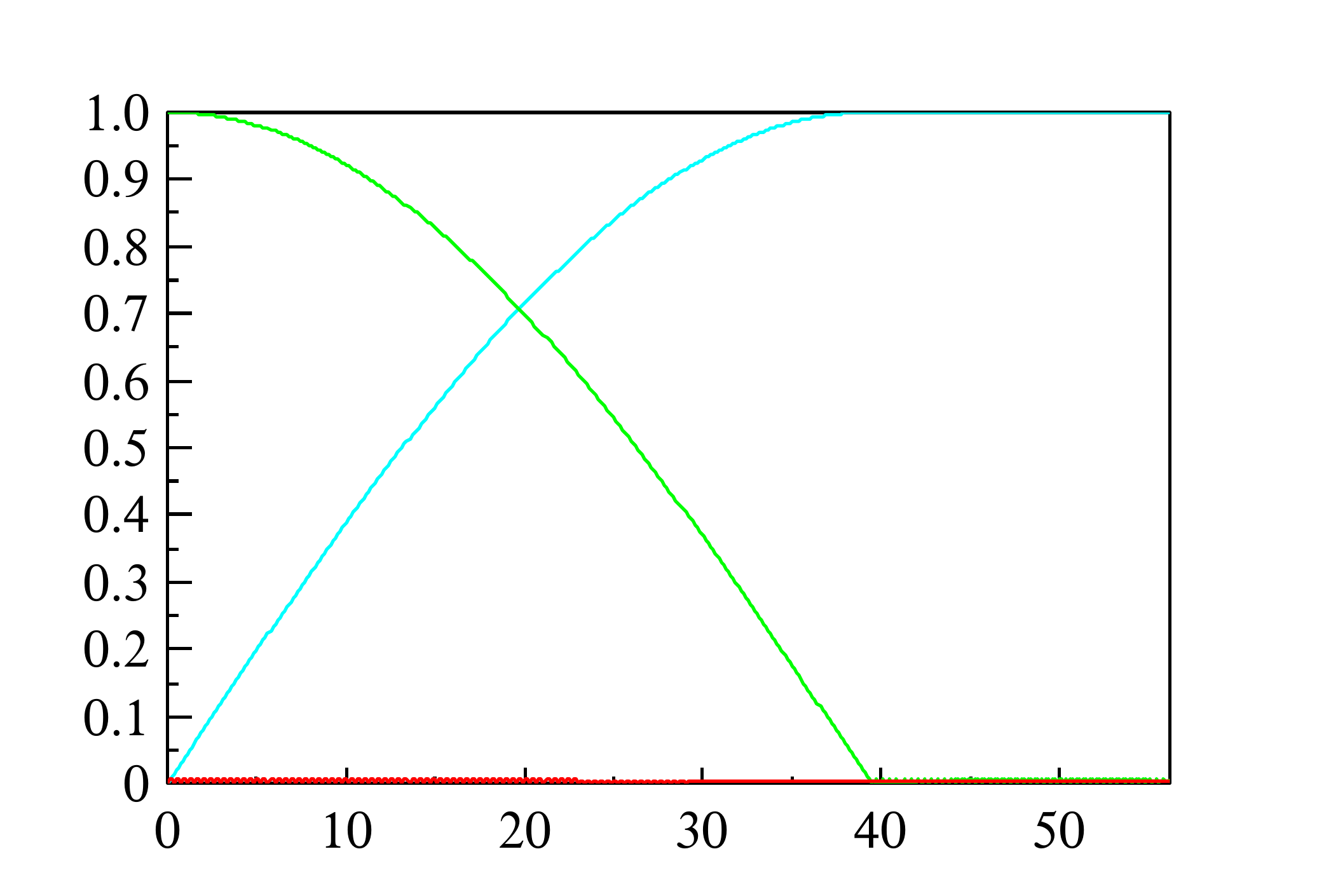}
 % HarmonicFromLevel1.pdf: 1276x696 pixel, 72dpi, 45.01x24.55 cm, bb=0 0 1276 696
 \caption{Time evolution of the moduli of the first three coordinates of $\Upsilon^u_t \phi_2$ in the case of the perturbed harmonic oscillator. First coordinate in blue, second coordinate in green, third coordinate in red. Time scale is 1:100.}
 \label{fig:2}
\end{figure}

\begin{figure}[tb]
 \centering
 \includegraphics[width=8.4cm]{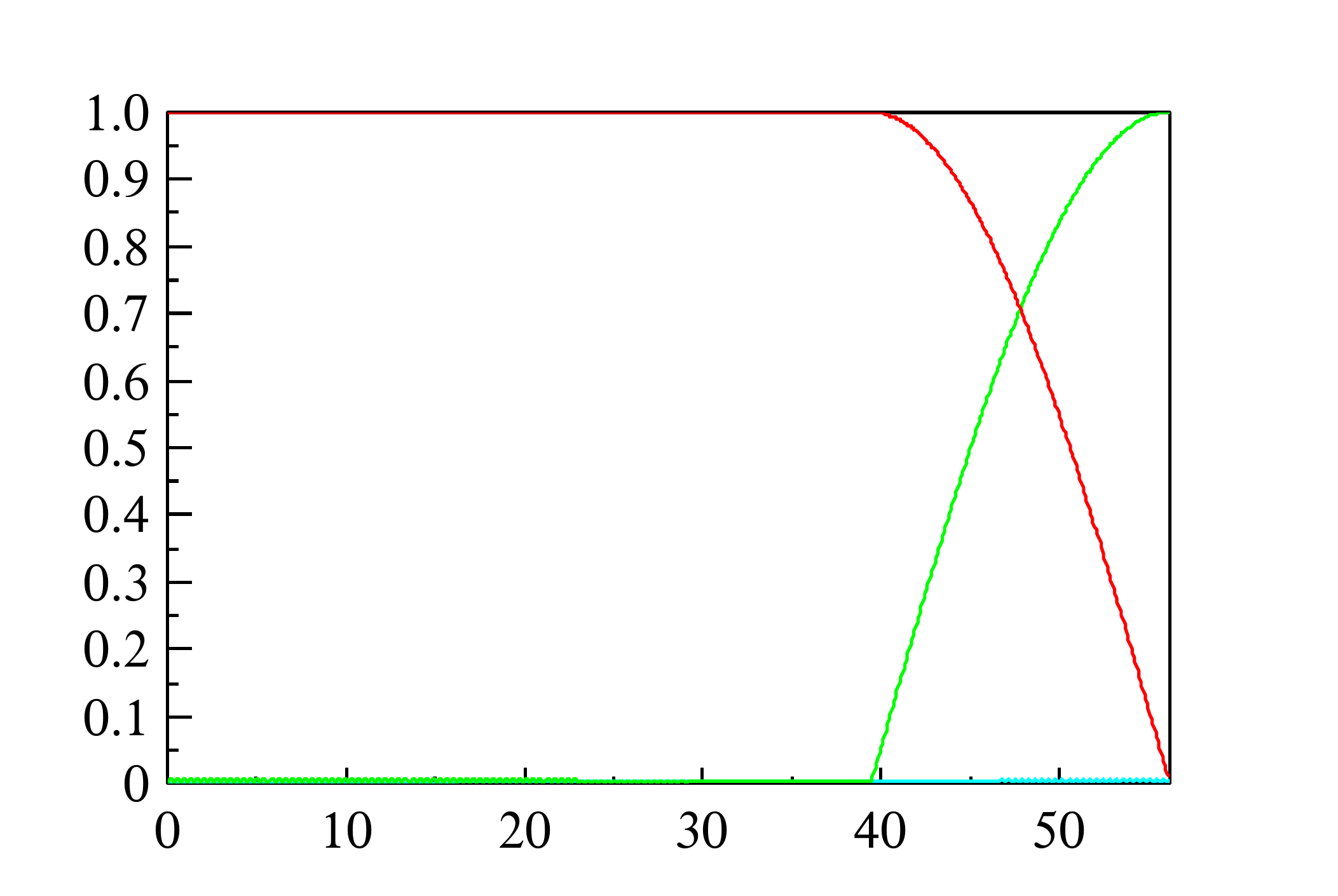}
 % HarmonicFromLevel1.pdf: 1276x696 pixel, 72dpi, 45.01x24.55 cm, bb=0 0 1276 696
 \caption{Time evolution of the moduli of the first three coordinates of $\Upsilon^u_t \phi_3$ in the case of the perturbed harmonic oscillator. First coordinate in blue, second coordinate in green, third coordinate in red. Time scale is 1:100.}
 \label{fig:3}
\end{figure}

\section{Particle in a box}\label{sec:potentialwell}

\subsection{Physical model}
We consider now the case of a particle confined in $(0,\pi)$. This model has been extensively studied
by several authors in the last few years and was the first quantum system for which a positive controllability result has been obtained. Beauchard proved exact controllability in some dense subsets of $L^2$ using Coron's return method (see \cite{beauchard-coron, beauchard-mirrahimi} for a precise statement). Nersesyan obtained approximate controllability results using Lyapunov techniques. In the following, we extend these controllability results to simultaneous controllability and provide some estimates of the $L^1$-norm of the controls achieving simultaneous controllability.% the transfer between two density matrices.

The Schr\"{o}dinger equation writes
\begin{equation}\label{EQ_potential_well}
 i \frac{\partial \psi}{\partial t}=-\frac{1}{2} \frac{\partial^2 \psi}{\partial x^2} - u(t) x \psi(x,t)
\end{equation}
with boundary conditions $\psi(0,t)=\psi(\pi,t)=0$ for every $t \in \mathbf{R}$.

In this case ${\cal H}=L^2 \left ( (0,\pi), \mathbb{C} \right )$ endowed with the Hermitian product $\langle \psi_1,\psi_2\rangle= \int_{0}^{\pi} \overline{\psi_1(x)} \psi_2(x) dx$. The operators $A$ and $B$ are defined by $A\psi= i\frac{1}{2} \frac{\partial^2 \psi}{\partial x^2}$ for every $\psi$ in
$D(A)= (H_2 \cap  H_0^1 )\left((0,\pi), \mathbb{C}\right)$, and $B\psi=i x \psi$.
An Hilbert basis of $H$ is $(\phi_k)_{k\in \mathbf{N}}$ with $\phi_k:x\mapsto \sin(kx)/\sqrt{2}$. For every $k$, $A\phi_k=-ik^2/2 \phi_k$.

For every $j,k$ in $\mathbf{N}$, %$j\neq k$,
$$
b_{jk}:=\langle \phi_j, B \phi_k\rangle =
%\left\{
\begin{cases}
(-1)^{j+k}\frac{ 2 jk}{(j^2-k^2)^2} & \mbox{ if } j-k \mbox{ odd}\\
0 & \mbox{otherwise}.
\end{cases}
%\right.
$$

Despite numerous degenerate transitions, the system is approximately simultaneously controllable (see \cite[Section 7]{Schrod2}).

\subsection{Estimates}

%In order to estimate the error done when replacing the original infinite dimensional system by its Galerkin approximation, we may be tempted to apply 
Using Proposition \ref{prop:gga} to estimate the 
error done when replacing infinite dimensional system by its Galerkin approximation
 one finds, for $\|u\|_{L^1}=9\pi/16$ (see \cite[Remark 4]{feps}, with $K=9\pi/16$, $d=\pi$, $k=1$, $r=1$, $c_1(A,B)\leq \pi+2$, $\varepsilon=10^{-3}$), that if $N> 1.6~10^7$, then 
$$
\|\pi_N \Upsilon^u_t\phi_1-X^u_{(N)}(t,0)\phi_1\|\leq 10^{-3}.
$$

This estimation is definitely too rough to allow easy numerical simulations: matrix $B^{(10^7)}$ has about $5~10^{13}$ non-zeros entries, the numerical simulations at such scale are difficult without large computing facilities.  We have to go more into details to obtain finer estimates.

Assume that, for some $N$ in $\mathbf{N}$ and $\eta>0$, the control $u:[0,T]\rightarrow \mathbf{R}$ is such that, for every $t$ in $[0,T]$,
$$\|X^{u}_{(N)}(t,0)\pi_3 -\pi_3X^{u}_{(N)}(t,0)  \|\leq \eta.$$ We have
\begin{eqnarray*}
\lefteqn{\pi_{3} X^{u}_{(N)}(t,s)-X^{u}_{(N)}(t,s)\pi_3}\\
 &= &\pi_{3}X^{u}_{(N)}(t,0)X^{u}_{(N)}(0,s)  -X^{u}_{(N)}(t,0)X^{u}_{(N)}(0,s)  \pi_{3} \\
&=&
X^{u}_{(N)}(t,0) (\pi_{3}X^{u}_{(N)}(0,s)-X^{u}_{(N)}(0,s)\pi_3) \\
&&\quad \quad+ (\pi_3X^{u}_{(N)}(t,0)-X^{u}_{(N)}(t,0)\pi_3)X^{u}_{(N)}(0,s)
\end{eqnarray*}
so that
\begin{equation}\label{EQ_major_crochet}
\|\pi_{3} X^{u}_{(N)}(t,s) - X^{u}_{(N)}(t,s) \pi_{3}\| \leq 2 \eta.
\end{equation}

Projecting \eqref{eq:variationofconstants} on the first $3$ components
we have, for $j=1,2,3$ that
\begin{eqnarray}
\lefteqn{\| \pi_{3} \Upsilon^{u}_{t}(\phi_{j}) - \pi_3 X^{u}_{(N)}(t,0)  \phi_{j} \|} \nonumber\\
&\leq&
 \int_{0}^{t} \| \pi_{3} X^u_{(N)}(t,s) \pi_N B(\mathrm{Id} - \pi_{N}) \Upsilon^{u}_{s}(\phi_{j})\|
u(s)  \mathrm{d}s \nonumber\\
&\leq & \int_{0}^{t} \|  X^u_{(N)}(t,s) \pi_3 B(\mathrm{Id} - \pi_{N}) \Upsilon^{u}_{s}(\phi_{j})\| u(s)  \mathrm{d}s \nonumber\\
&&{+ \int_{0}^{t}\!\!\! \| (\pi_{3} X^u_{(N)}(t,s)-X^u_{(N)}(t,s)\pi_3)\|  \|B\|
u(s)  \mathrm{d}s} \nonumber\\
& \leq & \left (\int_0^T |u(t)|\mathrm{d}t \right )  ( \|\pi_3 B (\mathrm{Id}-\pi_N)\|
\nonumber \\
&&\quad  + 2 \| B\| \sup_{t}\|\pi_{3} X^u_{(N)}(t,0)-X^u_{(N)}(t,0)\pi_3)\|). \label{EQ_majoration_mario}
\end{eqnarray}

By skew-adjointness, $\|\pi_3 B (\mathrm{Id}-\pi_N)\|=\|(\mathrm{Id}-\pi_N) B \pi_3 \|$.
This last quantity tends to zero, and we are able to give estimates of the convergence rate. Indeed,
\begin{eqnarray*}
\| (\mathrm{Id} -\pi_N) B  \phi_1 \|^2
& \leq & \sum_{k>N} \left|   \frac{2 k}{(k-1)^{2}(1+k)^{2}}\right|^{2}\\
& \leq & 4 \sum_{k>N} \frac{1}{(k-1)^6}\\
&\leq&   \frac{1}{(N-2)^5}.
\end{eqnarray*}
Similarly,
\begin{eqnarray*}
\| (\mathrm{Id} -\pi_N) B  \phi_2 \|^2 &\leq& \frac{\sqrt{2}}{(N-3)^5}\\
\|(\mathrm{Id} -\pi_N) B  \phi_3 \|^2 &\leq& \frac{2}{(N-4)^5}.
\end{eqnarray*}

The procedure to induce a given transformation, up to a given tolerance $\varepsilon>0$, on the space $\mathrm{span}\{\phi_1,\phi_2,\phi_3\}$ is the following:
\begin{enumerate}
 \item Use estimates given in \cite{Schrod2} to give an a priori upper bound $K$ on the $L^1$-norm of the controls one will use.
\item From $K$ and $\varepsilon$, find $N$ such that $K \|\pi_3 B (\mathrm{Id}-\pi_N)\|\leq \varepsilon/2$.
\item In the finite dimensional space $\mathrm{span}\{\phi_1,\ldots,\phi_N\}$, consider the bilinear system $\dot x=(A^{(N)}+uB^{(N)})x$ and find a control $u$ achieving the desired transition up to $\varepsilon/(2K)$ and such that $\|u\|_{L^1}\leq K$. This can be done using standard averaging procedures (see for instance \cite{Sanders}).
\item Use \eqref{EQ_majoration_mario} to get an upper bound of the distance of the trajectories of $(\Sigma_N)$ and the actual infinite dimensional system.
\end{enumerate}

\subsection{Numerical simulations}

We illustrate the above procedure on an example. Fix $\varepsilon=7~10^{-2}$. We would like to find $u:[0,T]\rightarrow \mathbf{R}$ such that $|\langle \phi_3,\Upsilon^u_T \phi_1 \rangle|>1-\varepsilon$, $|\langle \phi_1,\Upsilon^u_T \phi_2 \rangle|>1-\varepsilon$ and $|\langle \phi_2,\Upsilon^u_T \phi_3 \rangle|>1-\varepsilon$ at final time $T$. For this example, we are not interested in the respective phases but the method can easily be generalized to address this point (see Section~\ref{sec:improvements} below).

From \cite{Schrod2}, the transition can be achieved with controls of $L^1$-norm smaller than $5\pi/4(9/8 +25/24)$. Using controls with better efficiencies (as described in \cite{periodic}), we can use controls with $L^1$-norm smaller than
$2(9/8+25/4)=13/3$.

Using the above estimates, one sees that if $N=20$, then
 $$K \|\pi_3 B (\mathrm{Id}-\pi_N)\|\leq \frac{13}{3} \frac{\sqrt{2}}{(N-4)^{5/2}} \leq 6~10^{-3}.$$

Last, we define $u$ by $u(t)=\cos(3t)/20$ for $0\leq t \leq 72$ and $u(t)=\cos(5t)/20$ for $72 <t \leq T=138$. We check that $\int_0^T|u(t)|\mathrm{d}t \leq
13/3$.
One checks numerically that $\|\pi_3 X^u_{(20)}(t,0)-X^u_{(20)}(t,0)\pi_3 \|\leq 1.3~10^{-3}$ for $t\leq 138$. From \eqref{EQ_major_crochet}, we get, for every $t,s\leq T$
$$
\|\pi_3X^{u}_{(N)}(t,s)-X^{u}_{(N)}(t,s)\pi_3\|\leq 2.6~10^{-3}.$$
From \eqref{EQ_majoration_mario}, we have, for $j=1,2,3$,
\begin{eqnarray*}
\| \pi_{3} \Upsilon^{u}_{t}(\phi_{j}) - \pi_3 X^{u}_{(N)}(t,0)  \phi_{j} \|&\leq& \!\!\!\!\!\frac{13}{3}(6~10^{-3}+8.2~10^{-3})\\
&\leq &6.1~10^{-2}.
\end{eqnarray*}
Conclusion follows from the numerical computations
\begin{eqnarray*}
 |\langle \phi_3,X^u_{(20)}(T,0)\phi_1 \rangle|&\approx& 0.99924\\
 |\langle \phi_1,X^u_{(20)}(T,0)\phi_2 \rangle|&\approx& 0.99943\\
  |\langle \phi_2,X^u_{(20)}(T,0)\phi_3 \rangle|&\approx& 0.99949.
\end{eqnarray*}
The actual precision is likely much better than $6.1~10^{-2}$ which is known for sure. However, our estimates do not allow a better conclusion.

The evolutions with respect to the time of the moduli of the first coordinates of $X^u_{(20)} \phi_k$ for $k=1,2,3$ are represented in Figures~4, 5, and 6.

\begin{figure}[t!]
 \centering
 \includegraphics[width=8.4cm]{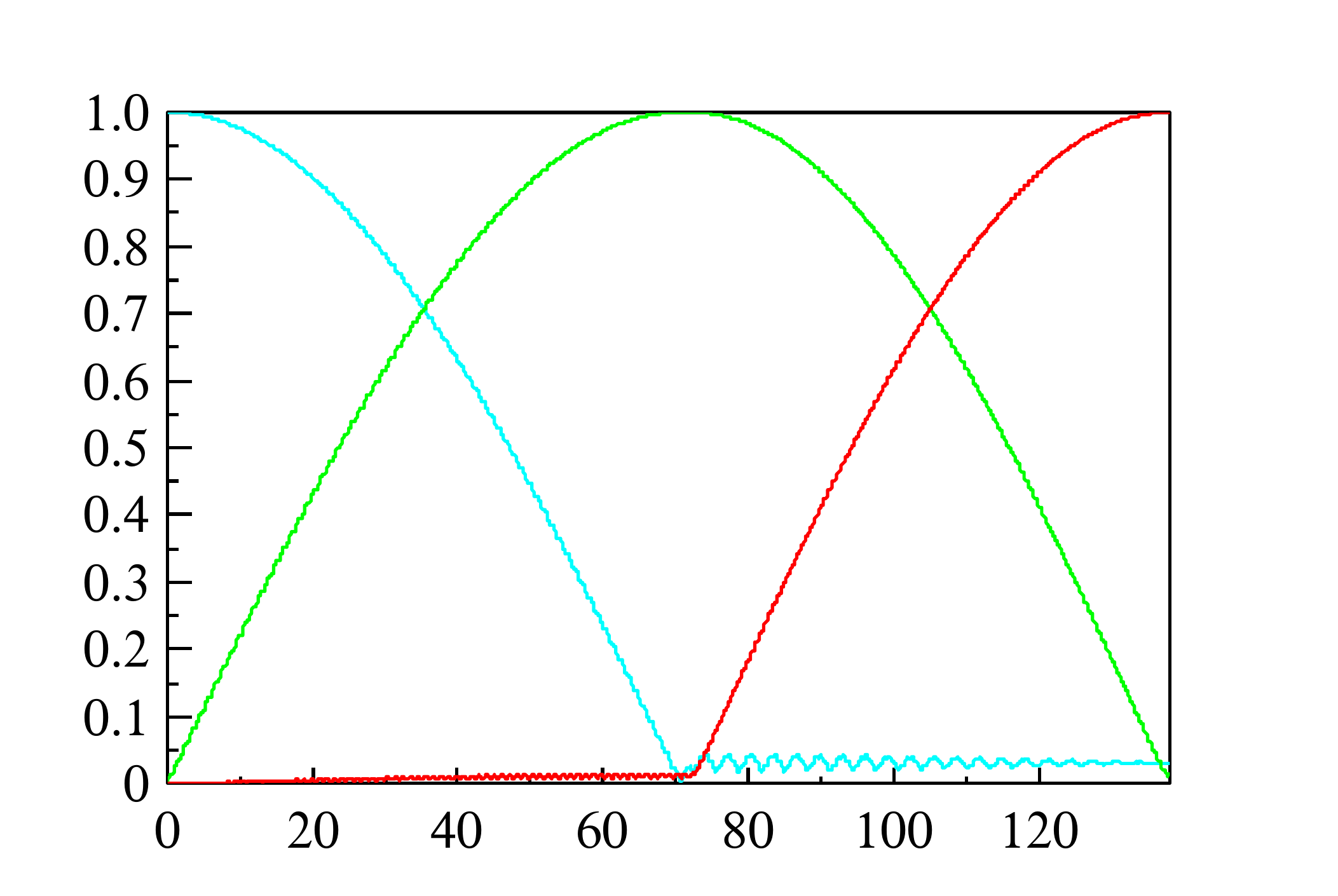}
  \caption{Time evolution of the moduli of the first three coordinates of $\Upsilon^u_t \phi_1$ in the case of the potential well. First coordinate in blue, second coordinate in green, third coordinate in red.}
 % FromLevel1.pdf: 1276x696 pixel, 72dpi, 45.01x24.55 cm, bb=0 0 1276 696
\end{figure}
\begin{figure}[t!]
 \centering
 \includegraphics[width=8.4cm]{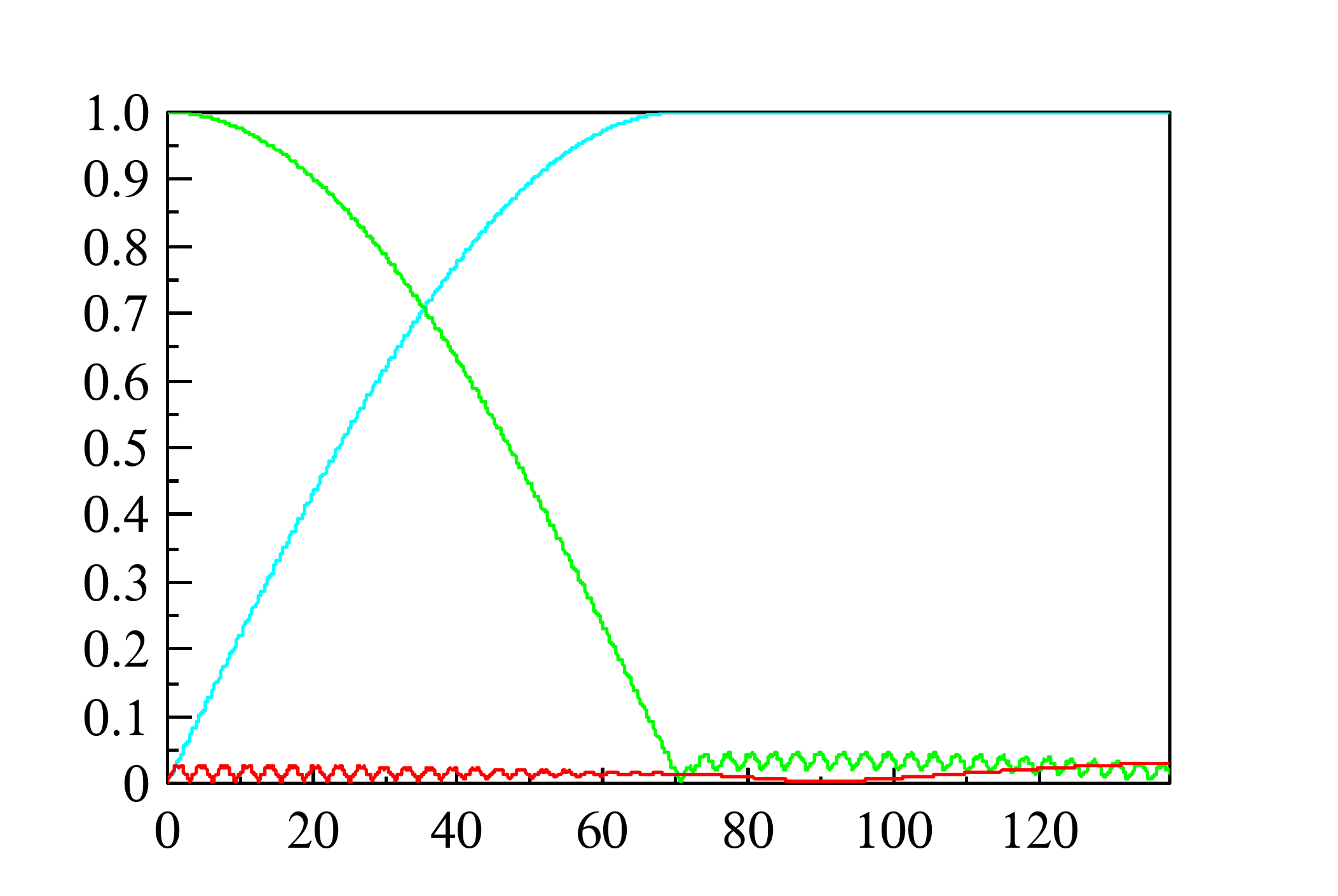}
\caption{Time evolution of the moduli of the first three coordinates of $\Upsilon^u_t \phi_2$ in the case of the potential well. First coordinate in blue, second coordinate in green, third coordinate in red.}
 % FromLevel1.pdf: 1276x696 pixel, 72dpi, 45.01x24.55 cm, bb=0 0 1276 696
\end{figure}
\begin{figure}[t]
 \centering
 \includegraphics[width=8.4cm]{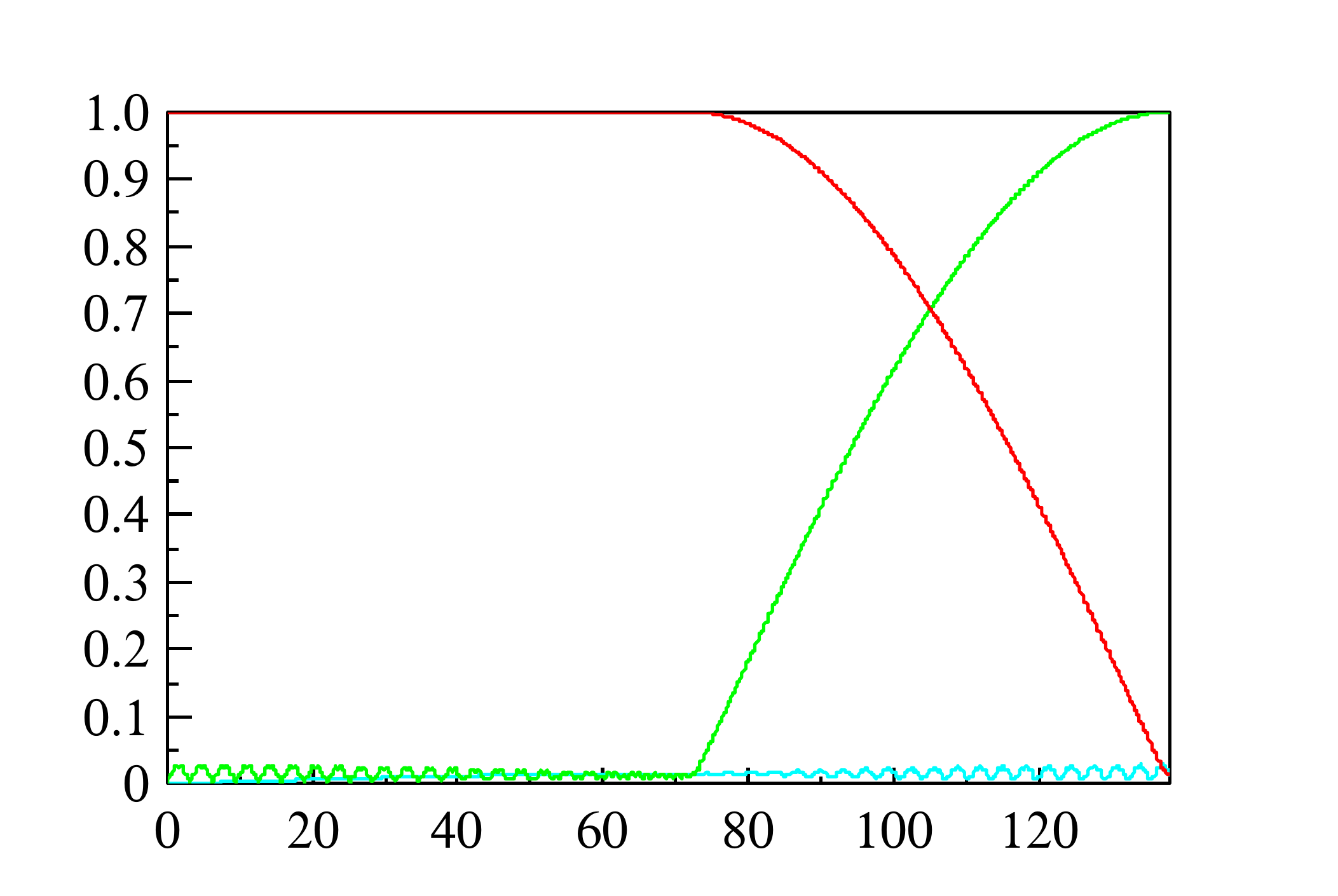}
\caption{Time evolution of the moduli of the first three coordinates of $\Upsilon^u_t \phi_3$ in the case of the potential well. First coordinate in blue, second coordinate in green, third coordinate in red.}
 % FromLevel1.pdf: 1276x696 pixel, 72dpi, 45.01x24.55 cm, bb=0 0 1276 696
\end{figure}

All the computations were done using the free software NSP, see \cite{NSP}. The source code for the simulation is available at \cite{source}. The total computation time is less than 4 minutes on a standard desktop computer.

\addtolength{\textheight}{-16.2cm}
\subsection{Possible improvements}\label{sec:improvements}

If one is interested not only in the modulus but also in the respective phases of the final points, it is enough to replace the functions $t\mapsto \cos(3t)/20$ and $t\mapsto \cos(5t)/20$ above by $t\mapsto \cos(3t+\theta_1)/20$ and $t\mapsto \cos(5t+\theta_2)/20$ respectively, where $\theta_1$ and $\theta_2$ are suitable phases.

In order to get better precision in the approximation  (i.e. a smaller $\varepsilon$), 
it is enough to replace the functions $t\mapsto \cos(3t)/20$ and $t\mapsto \cos(5t)/20$  above by
the functions $t\mapsto \cos(3t)/L$ and $t\mapsto \cos(5t)/L$ with $L$ large enough. The price to pay for a better precision is an increase in the time needed for the transfer.

\section{Conclusion and future works}
We have shown how it was possible to implement a quantum gate on two types of infinite dimensional quantum oscillators. Our method provides rigorous estimates and permits numerical simulations that can be run on standard desktop computers.

A limitation of our models is that the Schr\"{o}dinger equation neglects decoherence. This approximation may be justified for time small with respect to the relaxation time of the quantum system. Future works may focus on the optimization of the time of implementation.

\bibliographystyle{plain}
\bibliography{biblioACC}

\end{document}